\documentclass[12pt,oneside,reqno]{amsart}
\usepackage{amsmath,amssymb,amsfonts,latexsym,txfonts,pxfonts,wasysym,amsthm}
\usepackage{amssymb,longtable}
\usepackage{hyperref}
\usepackage{lscape,graphicx}%
\usepackage{amsmath}%
\usepackage{amsfonts}%
\usepackage{amssymb}%
\usepackage{graphicx}
\usepackage{times}
\usepackage{cite}
\allowdisplaybreaks \newmuskip\pFqmuskip
\newcommand*\pFq[6][8]{%
  \begingroup
  \pFqmuskip=#1mu\relax
    \mathcode`\.=\string"8000
   \begingroup\lccode`\~=`\.
  \lowercase{\endgroup\let~}\pFqdot
   {}_{#2}F_{#3}{\left[\genfrac..{0pt}{}{#4}{#5};#6\right]}%
  \endgroup
  }
\newcommand{\fs}{\frac{1}{2}}
\title[Comment on a paper ``Watson - like Formulae for terminating $_{3}F_2$ series" by Chu and Zhou]{Comment on a paper``Watson - like Formulae for terminating $_{3}F_2$- series"  by Chu and Zhou}
\author{\textbf{Arjun K. Rathie}}
\address{ Arjun K. Rathie, Department of Mathematics,  School of Mathematical and Physical Sciences, Central University of Kerala,  Tejaswini Hills, Periye P.O., Kasaragod, 671316,  Kerala State, India.}
\email { akrathie@cukerala.ac.in}
\begin{document}
\begin{abstract}
In a recent paper, Chu and Zhou [Advances in Combinatorics, I.S. Kotsireas and E.V. Zima(eds.), 139-159 (2013)] established in all 40 closed formulae   for terminating  Watson-like hypergeometric $_{3}F_2$- series by investigating through Gould and Hsu's fundamental pair of inverse series relations, the dual relations of Dougall's formula for the very well - poised $_{5}F_4$ - series.   

The aim of this short note is just to point out that out of 40 results, 33 results have already been discovered in 1992 by Lavoie, et al.
\\

\textbf{2010 Mathematics Subject Classifiction} :  33C20 

\textbf{Keywords} : Generalized Hypergeometric Functions, Watson Theorem, Dougall Theorem
\end{abstract}
\maketitle
\section{OBSERVATIONS}
\indent In a recent paper, Chu and Zhou[1] obtained 40 closed formulae for terminating Watson-like hypergeometric $_{3}F_2$ series in the form 
\begin{equation}
\pFq{3}{2}{-2n, \; a+2n, \; c\; }{\fs(a+i+1), \; 2c+j}{1} 
\end{equation}
and
\begin{equation}
\pFq{3}{2}{-2n-1, \; a+2n+1, \; c\; }{\fs(a+i+1), \; 2c+j}{1} 
\end{equation}
for $-5 \leq i,j \leq 5.$\\
by investigating through Gould and Hsu's fundametal pair of inverse series relations, the dual relations of Dougall's formula for the very well-poised $_{5}F_4$-series.\\

Where as, in 1992, Lavoie et al.[2] have already obtained explicit expressions of 

\begin{equation}
 \pFq{3}{2}{a, \; b, \; c\; }{\fs(a+b+i+1), \; 2c+j}{1} 
\end{equation}
for $i,j = 0, \pm 1, \pm 2$.\\
with the help of contiguous function relations for $_{3}F_2$.
In the same paper[2], they have also deduced explicit expressions of (1) and (2) for $i,j = 0, \pm1, \pm2$.

Therefore, it is observed here that out of 40 results obtained by chu and Zhou[1], only 7 results are new and rest 33 results are already recorded in [2].

Moreover, two results of the form
\begin{equation}
\pFq{3}{2}{a, \; b, \; c\; }{\fs(a+b+i+1), \; 2c}{1} 
\end{equation}
and \begin{equation}
\pFq{3}{2}{a, \; b, \; c\; }{\fs(a+b-i+1), \; 2c}{1} 
\end{equation}

in the most general form for any $i=0,1,2, \ldots$ are also recorded in [3].

\end{document}